\documentclass[conference]{IEEEtran}
\usepackage[T1]{fontenc}
\usepackage{stackengine}
\IEEEoverridecommandlockouts
\usepackage{cite}
\usepackage{amsmath,amssymb,amsfonts, array}
\usepackage{algorithmic}
\usepackage{graphicx}
\usepackage{textcomp}
\usepackage{xcolor}
\usepackage{algorithm,algorithmic}
\usepackage{mathtools, stackengine, mathrsfs, setspace, bm, multirow, tikz, pgfplots, adjustbox, calc, url, xtab}
\usepackage[colorlinks=true, citecolor=blue, linkcolor=black!20!red]{hyperref}
\usepackage[caption=false,font=footnotesize]{subfig}
\usepackage{xparse}
\usetikzlibrary{arrows,automata}

\allowdisplaybreaks

\begin{document}
\bstctlcite{IEEEexample:BSTcontrol}
\title{Risk-Averse Self-Scheduling of Storage in Decentralized Markets 
\thanks{Mr. Yurdakul gratefully acknowledges the support of the German Federal Ministry of Education and Research and the Software Campus program under Grant 01IS17052. Mr. Billimoria's Doctor of Philosophy at the University of Oxford is supported by the Commonwealth PhD Scholarship, awarded by the Commonwealth Scholarship Commission of the UK (CSCUK).}
}

\author{\IEEEauthorblockN{Ogun Yurdakul\IEEEauthorrefmark{1}\IEEEauthorrefmark{2} and Farhad Billimoria\IEEEauthorrefmark{3}
}
\IEEEauthorblockA{\IEEEauthorrefmark{1}Energy Systems and Infrastructure Analysis Division, Argonne National Laboratory, Lemont, IL 60439, USA.}
\IEEEauthorblockA{\IEEEauthorrefmark{2}Department of Electrical Engineering and Computer Science, Technical University of Berlin 10623 Berlin, Germany}
\IEEEauthorblockA{\IEEEauthorrefmark{3}Department of Engineering Science, University of Oxford, Oxford OX1 2JD, United Kingdom}
\vspace{-0.8cm}
}
\maketitle
\begin{abstract} 
Storage is expected to be a critical source of firming in low-carbon grids. A common concern raised from ex-post assessments is that storage resources can fail to respond to strong price signals during times of scarcity. While commonly attributed to forecast error or failures in operations, we posit that this behavior can be explained from the perspective of risk-averse scheduling. Using a stochastic self-scheduling framework and real-world data harvested from the Australian National Electricity Market, we demonstrate that risk-averse storage resources tend to have a myopic operational perspective, that is, they typically engage in near-term price arbitrage and chase only few extreme price spikes and troughs, thus remaining idle in several time periods with markedly high and low prices. This has important policy implications given the non-transparency of unit risk aversion and apparent risk in intervention decision-making.
\end{abstract}
\vspace{-0.3cm}
\section{Introduction}\label{sec1}
In renewable-dominated power systems, electricity storage is expected to play an important role in maintaining system balance and resource adequacy. In decentralized markets, system operators (SOs) do not issue commitment instructions, instead relying upon full-strength price formation and iterative pre-dispatch systems to signal the need for additional or reduced resource commitment.
%\footnote{For the avoidance of doubt however, the SO will still centrally clear a security-constrained economic dispatch.}
During periods of scarcity, understanding the projected available capacity of resources to service load at intraday system stress points is critical to the SO's decision-making on whether to intervene or apply emergency protocols, such as load shedding. \par
Storage resources introduce new complexity to the problem because while the resource may be operationally available, whether it can actually charge or discharge at a point in time depends upon its state of charge (driven by prior decisions) and technical constraints. As a case in point, the lack of transparency on the availability of energy-limited plants was one of the key drivers of the Australian Energy Market Operator's (AEMO) decision to suspend the National Electricity Market (NEM) in June 2022 during a period of extreme scarcity. In certain cases, storage resources were observed to not discharge during the highest prices of the day (while discharging at more muted prices earlier in the day).\footnote{Our online companion \cite{onlinecomp} provides an empirical record of many such observations of battery dispatch and energy prices for the NEM.} This seemingly non-intuitive behavior is commonly attributed to error or failure (e.g., in forecasting or operational management) of storage operations. Our work however proposes a second explanation: that of risk-averse self-scheduling under uncertainty, a concomitant of which is that such behavior could well be rational and utility-maximizing, and thus persist in energy markets. \par
In the literature, risk-averse self-scheduling of thermal resources has been extensively studied. Risk-constrained self-scheduling is explored in \cite{1318695,Papavasiliou2015} given the internalization of the non-convex attributes of thermal resources. In \cite{jabr2005robust}, the conditional Value-at-Risk measure is used to develop bidding curves derived from an optimal self-schedule. For storage and in relation to decentralized markets, inter-temporal tradeoffs are considered from the perspective of look-ahead bidding to optimize state-of-charge and maximize profits \cite{wang2017look}. In  \cite{kazempour2009risk} a risk-constrained problem is formulated for a pumped storage plant to trade-off between expected profit and risks in energy and ancillary service markets.\par
In this paper, we work out a novel risk-averse stochastic optimization framework for the self-scheduling of storage in decentralized markets under uncertain prices. In Section \ref{sec2}, we develop the mathematical background and market architecture, leading to the precise mathematical formulation of the risk-averse scheduling problem in Section \ref{sec3}. In Section \ref{sec4}, we conduct several experiments using real-world actual data from the NEM and storage resources bidding therein.\par
Using the insights we draw from the experiments, we lay out in Section \ref{sec4} two key contributions of this work. First, we provide a novel explanation for the seemingly non-intuitive behavior of storage. Specifically, we demonstrate how the increasing uncertainty of prices with longer horizons can lead a risk-averse decision-maker to adopt a rational but myopic approach to dispatch. We illustrate how such a decision-maker can favor low-risk near-term dispatch decisions at more moderate prices rather than higher-reward but more uncertain peaks and troughs. Second, we present valuable insights into the sensitivity of the expected profits to the duration and the degree of risk-aversion of a storage resource. We observe that while increasing the capacity of a storage resource can significantly boost profits for a risk-neutral decision-maker, it barely makes a dent for relatively more risk-averse decision-makers. We set out concluding remarks and policy implications in Section \ref{sec5}. 
\section{Mathematical Background}\label{sec2}
We consider a storage resource trading in a decentralized electricity market. The adopted market setting does not include a centrally organized short-term forward market for energy and relies solely on the real-time market (RTM) for the spot trading of energy, settled under a marginal pricing scheme\footnote{In this paper we are focused upon energy prices and do not consider, at this stage, ancillary service or other spot markets for which storage is eligible}. We assume that the resource acts as a price-taker, that is, it has no capability to exercise market power or alter market prices. \par
We denote by $k$ the index of each time period for which the RTM is cleared and by $K$ the total number of time periods in the problem horizon. We represent by $\kappa$ the duration of each time period $k$ in minutes, which is the smallest indecomposable unit of time in our analysis, during which we assume the system conditions hold constant. Using our notation, an RTM cleared for a day at half-hourly intervals would correspond to $\kappa=30$ min and $K=48$. Define the set $\mathscr{K} \coloneqq \{k \colon k = 1,\ldots,K\}$.\par
We denote by $\tilde{\lambda}_k, k \in \mathscr{K}$ as the uncertain energy price in time period $k$ at the transmission system bus or zone into which the storage resource is located. We construct the vector $\boldsymbol{\tilde{\lambda}}\coloneqq [\tilde{\lambda}_{1}\cdots\tilde{\lambda}_{k} \cdots \tilde{\lambda}_{K}]^{\mathsf{T}}$. We assume that the SO determines and publishes pre-dispatch energy prices over the $K$ time periods in order to provide market participants and itself with advance information necessary to plan the physical operation of the power system. We write the relation $\boldsymbol{\tilde{\lambda}} \coloneqq \boldsymbol{\bar{\lambda}} + \boldsymbol{\tilde{\epsilon}}$, where $\boldsymbol{\bar{\lambda}}$ denotes the vector of pre-dispatch energy prices and $\boldsymbol{\tilde{\epsilon}}$ is a random vector of forecast errors representing the difference between the pre-dispatch and market-clearing prices.\par 
Typically, the uncertain forecast error exhibits a greater degree of variance with extending horizon \cite{onlinecomp}. While the forecast error between pre-dispatch and market-clearing prices are likely small for the periods close to the dispatch interval, the pre-dispatch prices that eventuate may significantly deviate from the market price for forecast horizons well away from the actual dispatch time. We use a set of scenarios to model $\boldsymbol{\tilde{\epsilon}}$, where $\boldsymbol{{\epsilon}}^{\omega}$ denotes the vector of forecast errors in scenario $\omega \in \Omega$ and ${\pi}^{\omega}$ denotes its associated probability of realization. We define the vector of LMPs for each scenario by $\boldsymbol{{\lambda}}^{\omega} \coloneqq \boldsymbol{\bar{\lambda}} + \boldsymbol{{\epsilon}}^{\omega}$ and write $\boldsymbol{{\lambda}}^{\omega} \coloneqq [{\lambda}_{1}^{\omega}\cdots {\lambda}_{k}^{\omega} \cdots {\lambda}_{K}^{\omega}]$, where ${\lambda}_{k}^{\omega}$ denotes the energy price in time period $k$ of scenario $\omega$. \par
We next turn to the mathematical modeling of storage. We denote by $p^{c}_k$ and $p^{d}_k$ the charging and discharging power of the storage resource in time period $k$ with maximum values $p^{c}_{M}$ and $p^{d}_{M}$, respectively. To ensure the mutual exclusivity of the charging and discharging modes, and to reduce the computational burden of the problem, we draw upon the single binary variable storage formulation presented in \cite{ychen} and enforce the following constraints:
\begin{IEEEeqnarray}{l'l'l}
p^{d}_k - p_{M}^{d} u_k \leq 0,& p^{c}_k - p_{M}^{c} (1 - u_k) \leq 0 & \forall k \in \mathscr{K}.\label{plim}\\
u_k \in \{0, 1\}, & p^{c}_k,p^{c}_k\geq 0   & \forall k \in \mathscr{K}. \label{bin}
\end{IEEEeqnarray}
We denote the efficiency of the storage unit by $\eta$ (assuming symmetry between charging and discharging efficiencies). We represent by $E_k$ the stored energy level at the end of the time period $k$ and write the intertemporal operational constraints:
\begin{IEEEeqnarray}{l;l}
\hspace{-0.5cm}E_{k} = E_{k-1} -  \frac{1}{\eta} p_k^{d} \kappa \frac{1\text{ h}}{60\text{ min}}  + {\eta} p_k^{c}  \kappa \frac{1\text{ h}}{60\text{ min}} &\forall k \in \mathscr{K} \setminus 1,\label{intertempgen}\\
\hspace{-0.5cm}E_{k} = E_{o} -  \frac{1}{\eta} p_k^{d} \kappa \frac{1\text{ h}}{60\text{ min}}  + {\eta} p_k^{c}  \kappa \frac{1\text{ h}}{60\text{ min}} & \forall k \in \{1\},\label{intertempini}
\end{IEEEeqnarray}
where $E_{o}$ denotes the initial stored energy level at the beginning of the problem horizon. The stored energy level in each period $k$ is bounded from above and below by $E_{M}$ and  $E_{m}$, respectively:  
\begin{IEEEeqnarray}{c}
E_{m} \leq E_k \leq E_{M} \hspace{.3cm}\forall k \in \mathscr{K}.\label{enl}
\end{IEEEeqnarray}
\par
In the next section, we harness the mathematical models laid out in this section to work out a framework for the risk-averse self-scheduling problem of a storage resource.
\vspace{-0.2cm}
\section{Risk-Averse Self-Scheduling} \label{sec3}
\vspace{-0.1cm}
Absent a short-term forward market, the storage resource is directly exposed to the uncertainty in the RTM prices. A great deal of studies in the literature considers a risk-neutral market participant, for which the storage resource may bring about large losses in some scenarios as long as those are offset by even greater gains in other scenarios. Such studies however do not cater to risk-averse decisions-makers, who constitute the focus of this work, that prefer to ward off large losses independent of potential profits. \par
A widely used measure for incorporating risk into decision-making is the Value-at-Risk (VaR). For a specified risk confidence level $\alpha$ in $(0,1)$, $\alpha$-VaR is an upper estimate of losses that is exceeded with probability $1-\alpha$. We denote the $\alpha$-VaR of the loss associated with a decision $\boldsymbol{x}$ by $\zeta_{\alpha}(\boldsymbol{x})$. Despite presenting an intuitive representation of losses, VaR exhibits several undesirable properties, including not taking account of the losses suffered beyond $\zeta_{\alpha}(\boldsymbol{x})$, non-coherence, and non-convexity when computed using scenarios. Instead, we draw upon the CVaR measure to manage risk in our framework. \par
For continuous distributions, the $\alpha$-CVaR of the loss for a decision $\boldsymbol{x}$, which we denote by $\phi_{\alpha}(\boldsymbol{x})$, is the expected loss given that the loss is greater than or equal to $\zeta_{\alpha}(\boldsymbol{x})$. The definition of CVaR for discrete distributions is yet more subtle, which is the case in our framework as we rely on scenarios to represent uncertain prices. Rockafellar and Uryasev \cite{2002} define $\phi_{\alpha}(\boldsymbol{x})$ for general distributions as the weighted average of $\zeta_{\alpha}(\boldsymbol{x})$ and the expected loss strictly exceeding $\zeta_{\alpha}(\boldsymbol{x})$. \par
From a mathematical standpoint, CVaR presents several appealing features. Pflug \cite{pflug} shows that it is a coherent risk measure. Most notably, $\phi_{\alpha}(\boldsymbol{x})$ can be efficiently computed by minimizing a piecewise linear and convex function \cite[Theorem 10]{2002}, which can be cast as a linear programming (LP) problem by introducing an additional variable. \par
A key thrust of our framework is to hedge the decision-maker against the risk of incurring high charging costs. For an associated confidence level $\alpha$, the $\alpha$-CVaR of the uncertain charging cost over the problem horizon can be evaluated by solving the following LP problem:
%\sum_{k \in \mathscr{K}}
\begin{IEEEeqnarray}{l'l}
\hspace{0.0cm} \underset{z^{\omega}, \zeta}{\text{minimize}} & \mathcal{R}_{\alpha}(\boldsymbol{x}) \coloneqq \zeta+ \frac{1}{1-\alpha}\sum_{\omega \in \Omega} \pi^{\omega}z^{\omega} ,\label{cvaro}\\
\hspace{0.0cm}\text{subject to} & z^{\omega} \geq \sum_{k \in \mathscr{K}} {\lambda}^{\omega}_k p_k^{c} - \zeta,\; z^{\omega} \geq 0 \label{cvarc1}
\end{IEEEeqnarray}
where the decision vector $\boldsymbol{x}$ succinctly represents the storage variables $p_k^{c}$, $p_k^{d}$, and $E_k$ for all $k$ in $\mathscr{K}$, $\zeta$, and the auxiliary variable $z^{\omega}$ introduced for each $\omega$ in $\Omega$ and $\pi^{\omega}$ denotes the probability of the scenario $\omega \in \Omega$. \par
In conjunction with minimizing the risk of incurring high charging costs,\footnote{We would be remiss if we did not set forth that the resource faces a risk while discharging as well, that of recording low revenues if the RTM prices materialize at significantly lower levels than their pre-dispatch counterparts. We choose to omit this risk in this paper and defer it to our future work. 
%Given the page limit, incorporating said risk into this work would overcrowd the results, depriving us of the chance to conduct an in-depth analysis of charging risk and the sensitivity of the results to various factors examined in Section \ref{sec4}.
} the decision-maker seeks to maximize the expected profits of the resource over the $K$ time periods:
\begin{IEEEeqnarray}{l}
\mathcal{P}(\boldsymbol{x}) \coloneqq \sum_{\omega \in \Omega} \pi^{\omega} \Big[ \sum_{k=1}^{K} \lambda^{\omega}_k \big(p^{d}_k-p^{c}_k\big) \Big].\label{profit}
\end{IEEEeqnarray}
We adjust the trade-off between these seemingly conflicting objectives by minimizing the weighted combination of $\mathcal{R}(\boldsymbol{x})$ and $\mathcal{P}(\boldsymbol{x})$ subject to the storage constraints laid out in Section \ref{sec2} and the CVaR constraints presented in this section. The risk-averse self-scheduling (RASS) problem is expressed as:
\begin{IEEEeqnarray}{l'll}
\text{RASS}: & \text{minimize} \;&\hspace{0.5cm} -\mathcal{P}(\boldsymbol{x}) + \beta \mathcal{R}(\boldsymbol{x}),\nonumber\\
& \text{subject to} \;&\hspace{0.5cm}\eqref{plim}\text{--}\eqref{enl}, \eqref{cvarc1}. \nonumber
\end{IEEEeqnarray}
The weight parameter $\beta \in [0, \infty) $ tunes the decision-maker’s degree of risk-aversion. While increasing values of $\beta$ underscores the decision-maker’s desire to mitigate risk, driving down $\beta$ toward zero represents a more risk-neutral decision-maker. Setting $\beta=0$ implies that the sole objective of the decision-maker is to maximize her profits---independent of the risk that her decisions entail.\par
The RASS problem is solved on a rolling-window basis with a shrinking horizon. At the outset, it is solved for the time periods $k=1\cdots K$ before the uncertain price for any of the time periods is revealed, yet the optimal decisions for only the first time period are implemented. After the RTM price for $k=1$ is observed, the RASS problem is solved again for $k=2\cdots K$, this time implementing the optimal decision for only $k=2$, and so forth. The process repeats until the RASS problem is solved for $k=K$ with a single-period horizon. The binding decisions of each rolling window are passed along to the subsequent window by explicitly enforcing that the stored energy level at the end of the first, binding time period of each window be equal to the initial stored energy level at the beginning of the ensuing window.
\section{Case Study and Results}\label{sec4}
In this section, we conduct two case studies to gain insights into the self-scheduling of a storage resource under different pre-dispatch price signals. The price data of both case studies are from the NEM. As spelled out below, we use the actual pre-dispatch prices observed in the NEM in two representative days to form the vector of pre-dispatch prices in our experiments. To construct the scenarios $\boldsymbol{{\lambda}}^{\omega},\,\omega \in \Omega$, we draw upon the historical differences between the pre-dispatch and RTM prices that eventuated in the NEM across 2019, yielding 11,680 observations, from which we randomly select 100 observations to form the scenarios of each case study. The price data have a temporal granularity of 30 min, that is, $\kappa=30$ min and $K=48$. Commensurate with the price data, we harness the data of the storage resources currently bidding in the NEM \cite{AEMO2022} in our experiments. Unless otherwise stated, we pick the risk confidence level as $\alpha=0.95$ and the storage efficiency as $\eta=0.85$ in all experiments. The data and the source code of all experiments are furnished in the online companion \cite{onlinecomp}. 
\subsection{Case Study I}\label{sec4a}
\begin{figure}[h]
{\includegraphics[width= \linewidth]{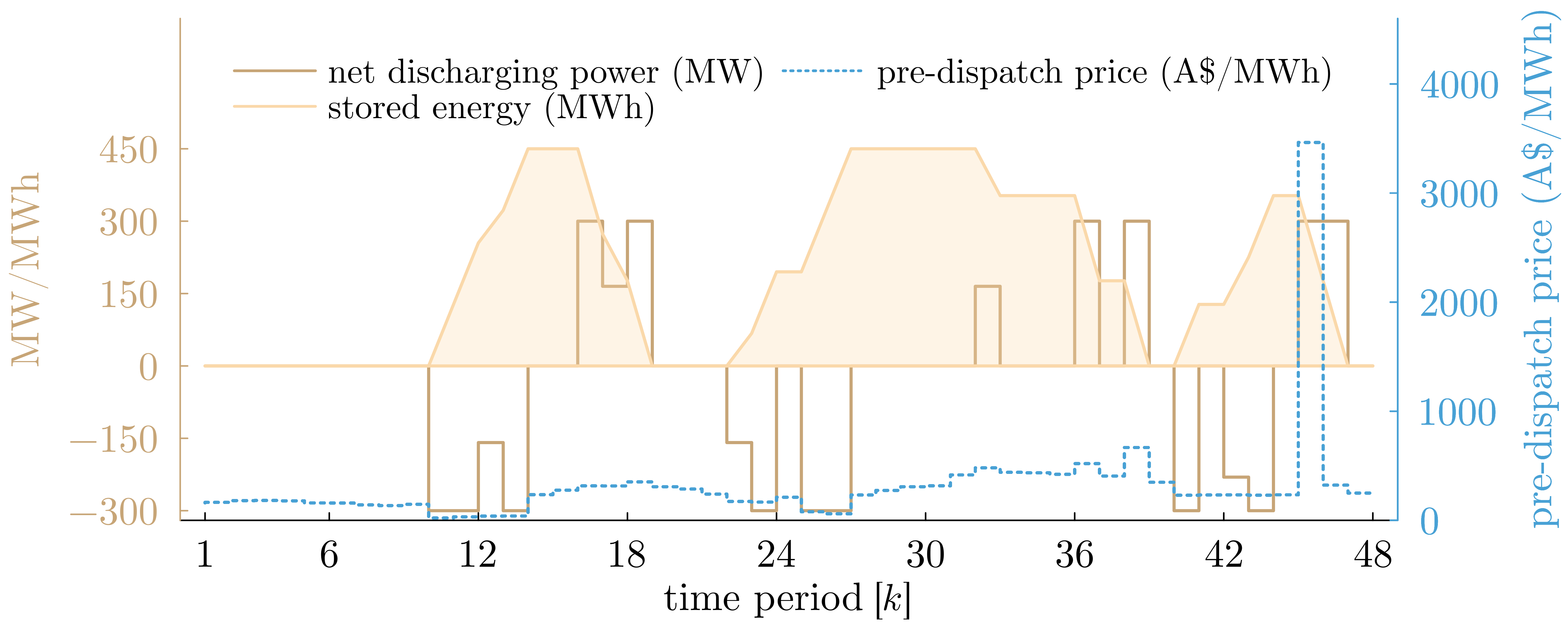}}\\
\vspace{-0.7cm}
\caption{Case Study I optimal storage dispatch decisions for $\beta=0$. Each time period $k$ has a duration of 30 minutes.}
\label{c1b0}
\vspace{-0.0cm}
\end{figure}

The storage data of Case Study I are from that of the Victorian Big Battery, which is a grid-connected storage resource in the Australian state of Victoria with a charging/discharging power limit of 300.0 MW and an energy storage capacity of 450.0 MWh. We use the pre-dispatch prices in Victoria for June 12, 2022, which constitutes the last day before the cumulative price threshold was exceeded in Victoria, triggering the onset of a series of market interventions that culminated in AEMO suspending the NEM on June 15, 2022. We start out by solving the RASS problem for $\beta = 0$ and $\beta=0.4$. We observe from Fig. \ref{c1b0} that the net discharging power (i.e., discharging power less charging power) under $\beta=0$, by and large, closely follows the pre-dispatch prices, effectively exploiting price differences across three time windows, with the first being between $k=11$ and $k=19$, 
%(that is, between 4:00 AM and 9:30 AM)
second between around $k=24$ and $k=36$, and finally that around the price spike at $k=46$.\par
\begin{figure}[h]
\vspace{-0.3cm}
{\includegraphics[width= \linewidth]{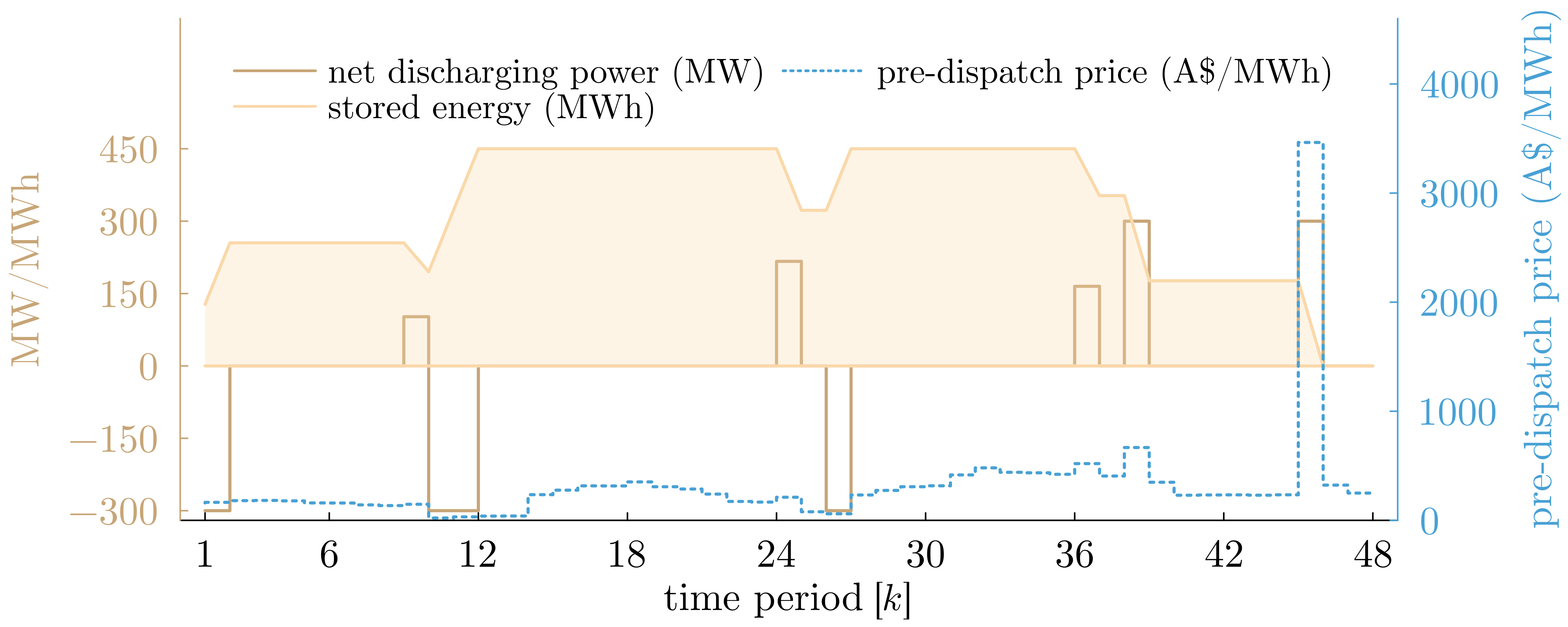}}\\
\vspace{-0.7cm}
\caption{Case Study I optimal storage dispatch decisions for $\beta=0.4$}
\label{c1b04}
\vspace{-0.1cm}
\end{figure}

We note from Fig. \ref{c1b04} that whereas the risk-neutral decision-maker ($\beta=0$) waits until the price reaches the early-hour minimum at $k=11$ to start charging, the risk-averse decision-maker ($\beta=0.4$) charges at its maximum limit right at the first two time periods, during which the pre-dispatch price is higher than that at $k=11$. We attribute this seemingly counter-intuitive behavior to the fact that when the RASS problem is solved at the beginning of the horizon, charging at $k=1$ entails a lower risk compared to $k=11$, as the market-clearing price for the initial time periods is expected hover closely around the pre-dispatch price. Since the uncertainty in forecast error increases with the length of the look-ahead period, the risk-averse decision-maker is driven to store more energy at the beginning of the horizon vis-\`a-vis the risk-neutral counterpart. \par
In this light, the risk-averse decision-maker incurs a greater \textit{lost opportunity cost} by storing a higher level of energy at $k=1$ and $k=2$ and discharging very little before $k=11$. This is because, by doing so, it reaches its maximum capacity at $k=12$ and foregoes its ability to store energy when the prices are around their lowest ebb of the day at $k=13$ and $k=14$. In contrast, the risk-neutral decision-maker manages to store 195 MWh of energy at $k=13$ and $k=14$. Tellingly, in order to store energy when the prices reach the minimum of the PM hours at $k=27$, the risk-averse decision maker winds up discharging at $k=25$. Although it could have well discharged between $k=15$ and $k=22$ at higher prices and could have recorded greater revenues, it prefers to discharge at a time period closer to that during which it charges again. Indeed, the dispatch decisions under $\beta=0.4$ is overall marred by a myopic focus. The risk averse decision-maker exploits price differences primarily at extremely short time windows, such as between $k=10$ and $k=11$ or between $k=25$ and $k=27$, unless the price significantly rises as around $k=46$. The prevailing myopic focus under $\beta=0.4$ can be ascribed to the increasing uncertainty in forecast error with longer horizons, driving the risk-averse decision-maker to arbitrage primarily between shorter windows. A ramification of this behavior is that the resource fails to respond to pre-dispatch price signals at time periods during which the system is in dire need. For instance, after the pre-dispatch price reaches A\$481.1/MWh at $k=33$ (more than a sevenfold increase from that at $k = 27$) the risk-averse decision-maker does not discharge any energy. Similarly, when the pre-dispatch price precipitously falls around $k=41$, the risk-neutral decision-maker stores 352.9 MWh, whereas the risk-averse decision-maker does not charge at all, failing to respond to the available 65.5\% price drop. As a result, after the price spikes at $k=46$, the risk-neutral decision-maker discharges 176.5 MWh more energy compared to the risk-averse decision-maker, greatly aiding the SO during a period of extreme scarcity.\par
\begin{figure}[h]
{\includegraphics[width= \linewidth]{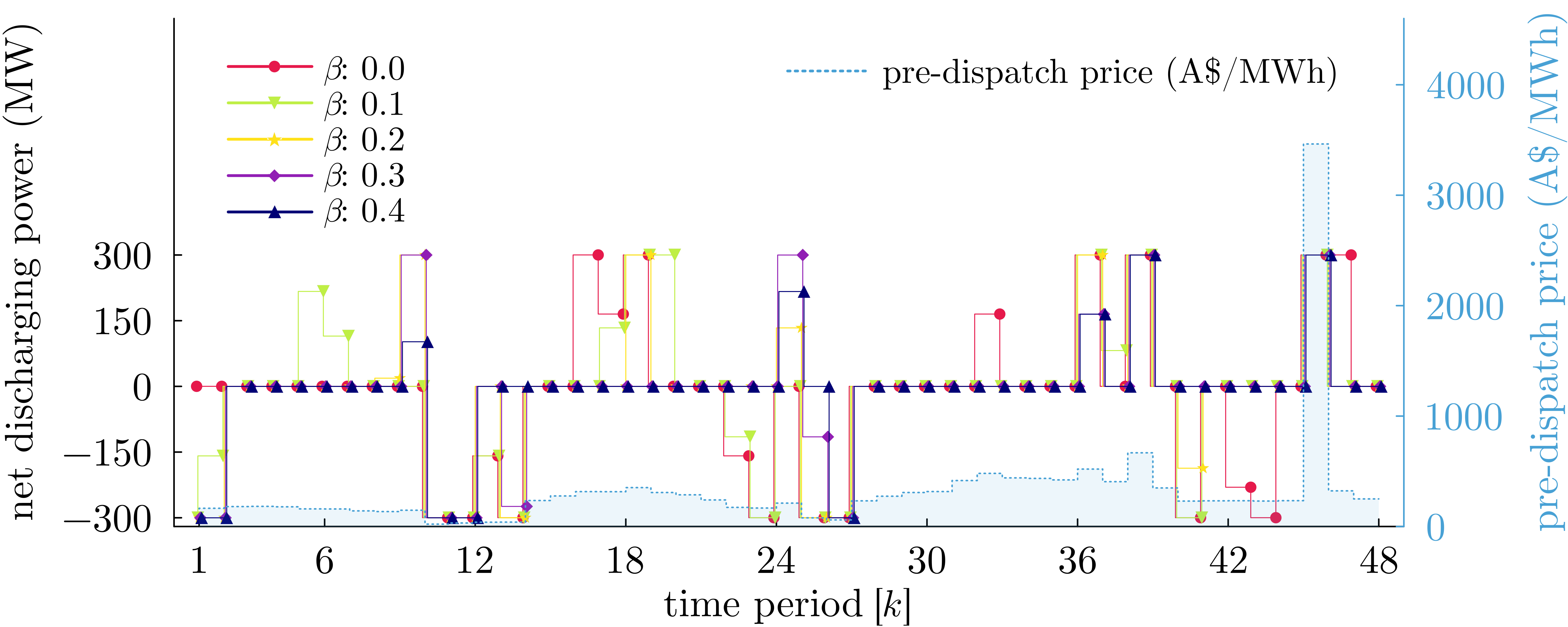}}\\
\vspace{-0.7cm}
\caption{Case Study I optimal net discharging levels (i.e., discharging less charging power) under different values of the weight parameter $\beta$.}
\label{c1a}
\vspace{-0.4cm}
\end{figure}
Fig. \ref{c1a} makes evident that the so-called myopic focus of the decision-maker becomes more conspicuous under increasing values of  $\beta$. Around the price spike at $k=19$, the storage resource discharges the highest level of energy under $\beta=0$, followed by $\beta=0.1$ and $\beta=0.2$, whereas no energy is discharged under $\beta=0.3$ and $\beta=0.4$. These decisions impart the relatively less conservative decision-makers ($\beta=0$, $\beta=0.1$, and $\beta=0.2$) the capability to store more energy during the price drop around $k=24$ compared to those under $\beta=0.3$ and $\beta=0.4$. These observations are echoed in the dispatch decisions throughout the rest of the day. For instance, as $\beta$ is reduced, the storage resource manages to more closely follow the pre-dispatch price signal before and around the sharp rise at $k=46$. Under decreasing values of $\beta$, the storage resource discharges a higher level of energy around $k=38$, gaining in turn the capability to store more energy during the dip in prices after $k=41$, which it can discharge when the pre-dispatch price abruptly climbs at $k=46$.
\subsection{Case Study II}\label{sec4b}
We next examine the battery dispatch decisions for a day in which the pre-dispatch prices are highly volatile. For this purpose, we draw upon the pre-dispatch prices in South Australia for January 16, 2019. The storage data are taken from that of the ESCRI storage resource, which is connected to the grid in South Australia and has a charging/discharging power limit of 30 MW and a capacity of 8 MWh.\par
\begin{figure}[h]
\vspace{-0.3cm}
{\includegraphics[width= \linewidth]{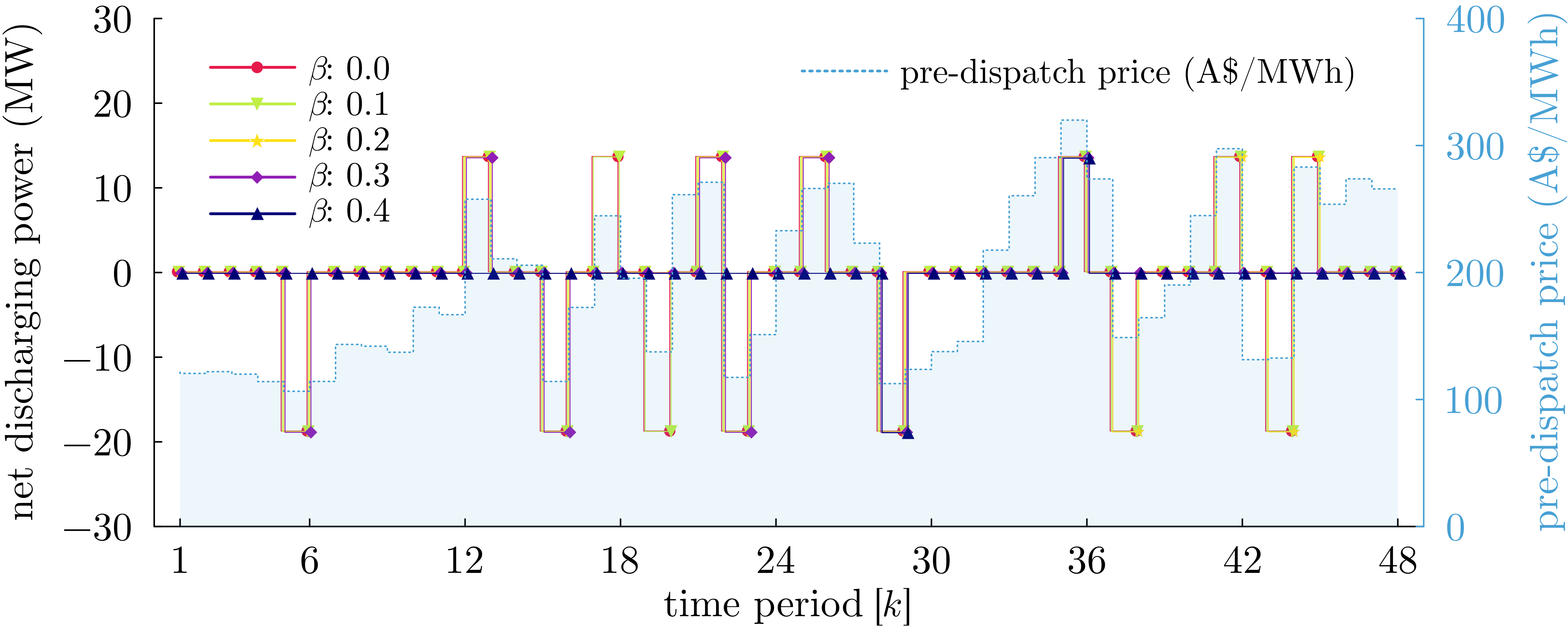}}\\
\vspace{-0.4cm}
\caption{Case Study II optimal net discharging levels under different values of the weight parameter $\beta$ for $\alpha=0.95$.}
\vspace{-0.2cm}
\label{c2a}
\end{figure}
We begin by solving the RASS problem by increasing $\beta$ from $0$ to $0.5$ in $0.1$ increments. The results in Fig. \ref{c2a} show that the relatively more risk-neutral cases ($\beta=0$ and $\beta=0.1$) closely follow the pre-dispatch price signals, charging/discharging at all of the seven price troughs/peaks of the day. Under $\beta=0.2$ and $\beta=0.3$, however, the storage resource prefers to arbitrage between only six and four time windows, respectively, whereas under $\beta=0.4$, it capitalizes on only the widest price spread, charging and discharging once throughout a highly volatile day. Most notably, under $\beta=0.4$, the resource relinquishes the chance to take advantage of a $56.7\%$ price difference between $k=22$ and $k=23$ when the price falls precipitously, yet changes course three periods later and records a $126.7\%$ rise between $k=23$ and $k=26$.\par
% As such, less risk-averse policies gain the chance to notch greater profits, which jibes with the results in Fig. \ref{c2ha}.\par
\begin{figure}[h!]
\centering
{\includegraphics[width= \linewidth]{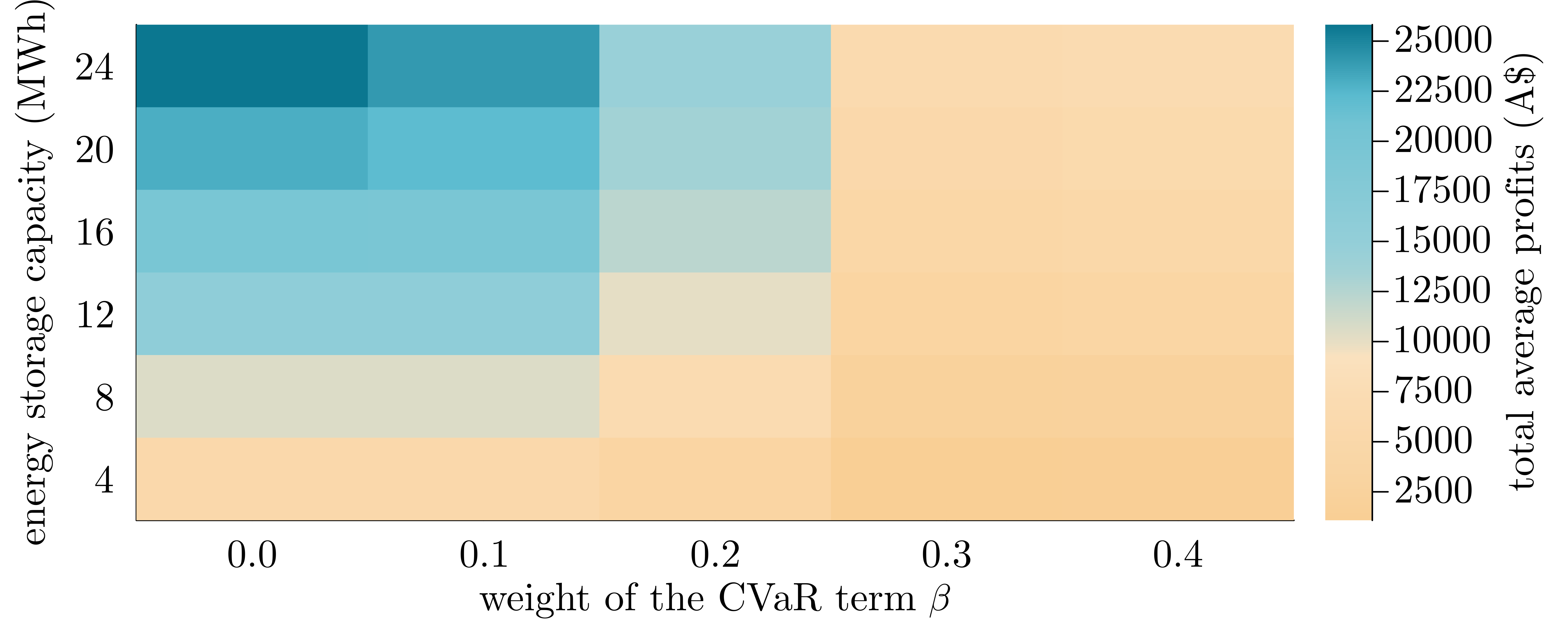}}\\
\vspace{-0.2cm}
\caption{Total expected profits under different values of $\beta$ and $E_M$.}
\label{c2hs}
\vspace{-0.2cm}
\end{figure}
We observe that the storage fails to reach its maximum charging and discharging power in the simulations, because if it were to sustain its maximum charging power over 30 minutes, it would exceed its total energy capacity. As such, we explore how the recorded profits would have evolved had the storage resource had a different capacity. To this end, we repeat the experiments by varying $E_M$ from $4.0$ MWh to $24.0$ MWh with $4.0$ MWh increments. We note from Fig. \ref{c2hs} that, across most values of $\beta$, the total expected profits rise under growing values of $E_M$, which can be ascribed to an increased ability to leverage price differences under a larger storage capacity. As $\beta$ increases, however, the profits seem to be plagued by diminishing returns, and indeed plateau at a certain $E_M$ level as $\beta$ increases beyond $0.2$. These observations bring home that risk attitudes (via increasing $\beta$ values) can create higher hurdles to notching greater profits, which are not necessarily mitigated by a larger storage capacity.\par
\begin{figure}[h!]
\centering
{\includegraphics[width= \linewidth]{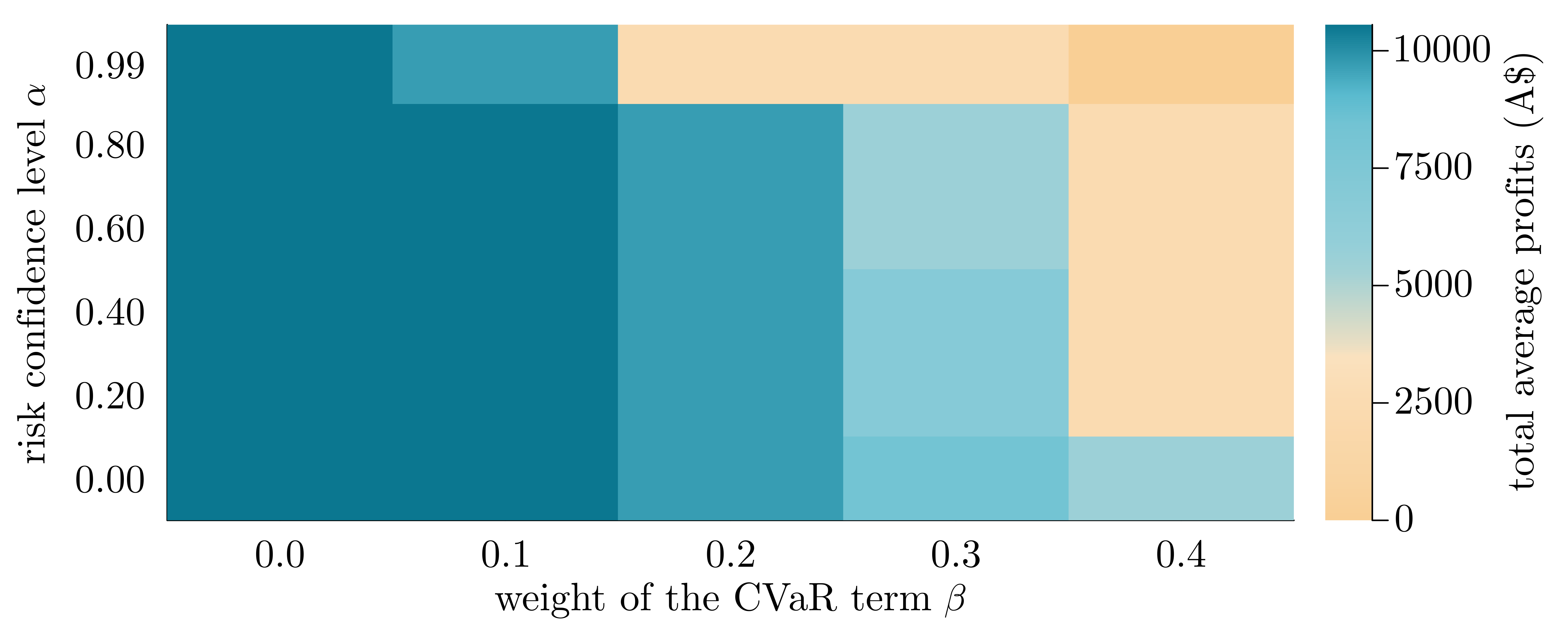}}\\
\vspace{-0.2cm}
\caption{Total expected profits under different values of $\alpha$ and $\beta$.}
\label{c2h}
\vspace{-0.2cm}
\end{figure}
In Fig. \ref{c2h}, we examine the influence of the risk confidence level $\alpha$ as well as $\beta$ on the total expected profits over the $K$ time periods. Note that, driving down $\alpha$ toward zero brings $\mathcal{R}_{\alpha}(\boldsymbol{x})$ toward the expected value of the charging costs, signifying a more risk-neutral decision-maker. In contrast, increasing $\alpha$ toward one drives $\mathcal{R}_{\alpha}(\boldsymbol{x})$ towards the highest value the costs can take,
%, i.e., $\underset{\omega \in \Omega}{\max}\big\{\sum_{k \in \mathscr{K}}\lambda_{k}^{c}p_{k}^{c}\big\}$, 
representing a more conservative decision-maker. Fig. \ref{c2h} bears out the diminished capability of risk-averse decision-maker to respond to pre-dispatch price signals and to arbitrage, manifesting itself through dwindling profits under increasing values of $\beta$ and $\alpha$. These results reaffirm the well-known relationship between risk and profit, whereby expected profits climb under less conservative decisions, characterized by the values of $\alpha$ and $\beta$ approaching zero.
\section{Conclusion and Policy Implications}\label{sec5}
\vspace{-.05cm}
We model the self-scheduling problem of a risk-averse  storage resource in decentralized markets. Four core results are gleaned from the analysis.  First that risk-aversion tends to lead to a myopic approach to dispatch most notably evident in seeking to arbitrage between near-term intervals, rather than taking charging positions in the expectation of more profitable but risky revenues much later in the look-ahead period.  Second that risk aversion tends to lead to reduced profits because of the conservative operating stance adopted by the resource.  Third, that this conservatism can mean that the resource forgoes opportunities to dispatch at peak spot prices, which are times of highest system scarcity. Finally, and importantly, while higher energy storage capacity can mitigate low profits, increasing the storage capacity has virtually no impact on profitability for higher risk-aversion levels.

There are important policy implications that flow from this study. First, SOs need to be acutely aware of the role of risk aversion in storage resource dispatch and bidding. In particular, higher levels of risk aversion may mean that storage resources may not be available for dispatch at times where scarcity price signals are most evident. As risk aversion itself is a non-observable quantity, this introduces uncertainty into an SO's short-term and medium term reliability forecasts, and risk into its decision-making on market intervention. This points to an increasing need for system transparency on key storage parameters including state-of-charge so that SOs have a better view on available storage capacity. Finally, it points to the need for the SO itself to develop a set of tools that quantify such risks and allow for better risk-adjusted decisions on intervention and other system operations during scarcity. 
\vspace{-.05cm}
\bibliographystyle{IEEEtran}
\bibliography{IEEEabrv,selfcomm}
\end{document}